\documentclass[]{article}
\usepackage{amsmath}
\usepackage{amssymb}
\usepackage{verbatim}
\usepackage{natbib}
\usepackage{fullpage}

\title{Fast Convergence of Stochastic Gradient Descent under a Strong Growth Condition} 
\author{Mark Schmidt and Nicolas Le Roux}

\begin{document}

\maketitle

\abstract{We consider optimizing a function smooth convex function $f$ that is the average of a set of differentiable functions $f_i$,
under the assumption considered by~\citet{solodov1998incremental} and~\citet{tseng1998incremental} that the norm of each gradient $f_i'$ is bounded
by a linear function of the norm of the average gradient $f'$. We show that under these assumptions the basic stochastic gradient method with 
a sufficiently-small constant step-size has an $O(1/k)$ convergence rate, and has a linear convergence rate if $g$ is strongly-convex.}

\section{Deterministic vs. Stochastic Gradient Descent}

We consider optimizing a function $f$ that is the average of a set of differentiable functions $f_i$,
\begin{equation}
\label{eq:one}
\min_{x \in \mathbb{R}^P} f(x) := \frac{1}{N}\sum_{i=1}^N f_i(x),
\end{equation}
where we assume that $f$ is convex and its gradient $f'$ is Lipschitz-continuous with constant $L$, meaning that for all $x$ and $y$ we have
\[
||f'(x) - f'(y)|| \leq L||x-y||.
\]
If $f$ is twice-differentiable, these assumptions are equivalent to assuming that the eigenvalues of the Hessian $f''(x)$ are bounded between
$0$ and $L$ for all $x$.

\emph{Deterministic gradient} methods for problems of this form use the iteration
\begin{equation}
\label{eq:four}
x_{k+1} = x_k - \alpha_k f'(x_k),
\end{equation}
for a sequence of step sizes $\alpha_k$.  In contrast,~\emph{stochastic gradient} methods use the iteration
\begin{equation}
x_{k+1} = x_k - \alpha_k f_i'(x_k),
\label{eq:three}
\end{equation}
for an individual data sample $i$ selected uniformly at random from the set $\{1,2,\dots,N\}$.  

The stochastic gradient method is appealing because the cost of its iterations is \emph{independent of} $N$.  However, in order to guarantee convergence stochastic gradient methods require a decreasing sequence of step sizes $\{\alpha_k\}$ and this leads to a slower convergence rate.  In particular, for convex objective functions the stochastic gradient method with a decreasing sequence of step sizes has an expected error on iteration $k$ of $O(1/\sqrt{k})$~\cite[\S14.1]{nemirovski1994efficient}, meaning that
\[
\mathbb{E}[f(x_k)] - f(x^*) = O(1/\sqrt{k}).
\]
In contrast, the deterministic gradient method with a \emph{constant} step size has a smaller error of $O(1/k)$~\cite[\S2.1.5]{nesterov2004introductory}.  The situation is more dramatic when $f$ is \emph{strongly} convex, meaning that
\begin{equation}
\label{eq:twelve}
f(y) \geq f(x) + \langle f'(x),y-x\rangle + \frac{\mu}{2}||y-x||^2,
\end{equation}
for all $x$ and $y$ and some $\mu > 0$.  For twice-differentiable functions, this is equivalent to assuming that the eigenvalues of the Hessian are bounded below by $\mu$.  For strongly convex objective functions, the stochastic gradient method with a decreasing sequence of step sizes has an error of $O(1/k)$~\cite[\S2.1]{nemirovski2009robust} while the deterministic method with a constant step size has an \emph{linear} convergence rate.  In particular, the deterministic method satisfies
\[
f(x_k) - f(x^*) \leq \rho^k[f(x_0) - f(x^*)],
\]
for some $\rho < 1$~\cite[\S8.6]{luenberger2008linear}.

The purpose of this note is to show that, if the individual gradients $f_i'(x_k)$ satisfy a certain strong growth condition relative to the full gradient $f'(x_k)$, the stochastic gradient method with a sufficiently small constant step size achieves (in expectation) the convergence rates stated above for the deterministic gradient method.

\section{A Strong Growth Condition}

The particular condition we consider in this work is that for all $x$ we have
\begin{equation}
\label{eq:two}
\max_i \{||f_i'(x)||\} \leq B||f'(x)||,
\end{equation}
for some constant $B$.  This condition states that the norms of the gradients of the individual functions are bounded by a linear function of the norm of the average gradient.  Note that this condition is very strong and is not satisfied in most applications.  In particular, this condition requires that any optimal solution for problem~\eqref{eq:one} must also be a stationary point for each $f_i(x)$, so that
\[
(f'(x) = 0) \Rightarrow (f_i'(x) = 0), \forall_i.
\]
In the context of non-linear least squares problems this condition requires that all residuals be zero at the solution, 
a property that can be used to show local superlinear convergence of Gauss-Newton algorithms~\cite[\S1.5.1]{bertsekas1999nonlinear}.

Under condition~\eqref{eq:two},~\citet{solodov1998incremental} and~\citet{tseng1998incremental} have analyzed convergence properties of \emph{deterministic incremental gradient} methods.  In these methods, the iteration~\eqref{eq:three} is used but the data sample $i$ is chosen in a deterministic fashion by proceeding through the samples in a cyclic order.  Normally, the deterministic incremental gradient method requires a decreasing sequence of step sizes to achieve convergence, but Solodov shows that under condition~\eqref{eq:two} the deterministic incremental gradient method converges with a sufficiently small constant step size.  Further, Tseng shows that a deterministic incremental gradient method with a sufficiently small step size may have a form of linear convergence under condition~\eqref{eq:two}.  However, this form of linear convergence treats full passes through the data as iterations, similar to the deterministic gradient method. 
 Below, we show that the stochastic gradient descent method achieves a linear convergence rate in expectation, using iterations that only look at one training example.

\section{Error Properties}
\label{sec:3}

It will be convenient to re-write the stochastic gradient iteration~\eqref{eq:three} in the form
\begin{equation}
\label{eq:seven}
x_{k+1} = x_k - \alpha(f'(x_k) + e_k),
\end{equation}
where we have assumed a constant step size $\alpha$ and where the error $e_k$ is given by
\begin{equation}
\label{eq:five}
e_k = f_i'(x_k) - f'(x_k).
\end{equation}
That is, we treat the stochastic gradient descent iteration as a full gradient iteration of the form~\eqref{eq:four} but with an error $e_k$ in the gradient calculation.  Because $i$ is sampled uniformly from the set $\{1,2,\dots,N\}$, note that we have
\begin{equation}
\label{eq:six}
\mathbb{E}[f_i'(x_k)] = \frac{1}{N}\sum_{i=1}^Nf_i'(x_k) = f'(x_k),
\end{equation}
and subsequently that the error has a mean of zero,
\begin{equation}
\label{eq:eight}
\mathbb{E}[e_k] = \mathbb{E}[f_i'(x_k) - f'(x_k)] = \mathbb{E}[f_i'(x_k)] - f'(x_k) = 0.
\end{equation}
In addition to this simple property, our analysis will also use a bound on the variance term $\mathbb{E}[||e_k||^2]$ in terms of $||f'(x_k)||$.  To obtain this we first use~\eqref{eq:five}, then expand and use~\eqref{eq:six}, and finally use our assumption~\eqref{eq:two} to get
\begin{equation}
\label{eq:nine}
\begin{aligned}
\mathbb{E}[||e_k||^2] & = \mathbb{E}[||f_i'(x_k) - f'(x_k)||^2]\\
& = \mathbb{E}[||f_i'(x_k)||^2 - 2\langle f_i'(x_k),f'(x_k)\rangle + ||f'(x_k)||^2]\\
& = \mathbb{E}[||f_i'(x_k)||^2] - 2\langle \mathbb{E}[f_i'(x_k)],f'(x_k)\rangle + ||f'(x_k)||^2\\
& = \frac{1}{N}\sum_{i=1}^N[||f_i'(x_k)||^2] - ||f'(x_k)||^2\\
& \leq (B^2 - 1)||f'(x_k)||^2.
\end{aligned}
\end{equation}

\section{Upper Bound on Progress}

We first review a basic inequality for inexact gradient methods of the form~\eqref{eq:seven}, when applied to functions $f$ that have a Lipschitz continuous gradient.  In particular, because $f'$ is Lipschitz-continuous, we have for all $x$ and $y$ that
\[
f(y) \leq f(x) + \langle f'(x),y-x\rangle + \frac{L}{2}||y-x||^2.
\]
Plugging in $x = x_k$ and $y = x_{k+1}$ we get
\[
f(x_{k+1}) \leq f(x_k) + \langle f'(x_k),x_{k+1}-x_k\rangle + \frac{L}{2}||x_{k+1}-x_k||^2.
\]
From~\eqref{eq:seven} we have that $(x_{k+1} - x_k) = -\alpha(f'(x_k) + e_k)$, so we obtain
\begin{equation}
\label{eq:13}
\begin{aligned}
f(x_{k+1}) & \leq f(x_k) - \alpha\langle f'(x_k),f'(x_k) + e_k\rangle + \frac{\alpha^2L}{2}||f'(x_k) + e_k||^2\\
& = f(x_k) - \alpha(1 - \frac{\alpha L}{2})||f'(x_k)||^2 - \alpha(1 - \alpha L)\langle f'(x_k),e_k\rangle + \frac{\alpha^2L}{2}||e_k||^2.
\end{aligned}
\end{equation}

\section{Descent Property}

We now show that, if the step size $\alpha$ is sufficiently small and the error is as described in Section~\ref{sec:3}, the expected value of $f(x_{k+1})$ is less than $f(x_k)$.  
In particular, we take the expectation of both sides of~\eqref{eq:13} with respect to $e_k$, and use~\eqref{eq:eight} and~\eqref{eq:nine} to obtain
\begin{equation}
\label{eq:ten}
\begin{aligned}
\mathbb{E}[f(x_{k+1})] & \leq f(x_k) - \alpha(1-\frac{\alpha L}{2})||f'(x_k)||^2 - \alpha(1-\alpha L)\langle f'(x_k),\mathbb{E}[e_k]\rangle + \frac{\alpha^2L}{2}\mathbb{E}[||e_k||^2]\\
& \leq f(x_k) - \alpha(1-\frac{\alpha L}{2})||f'(x_k)||^2  + \frac{\alpha^2L(B^2-1)}{2}||f'(x_k)||^2\\
& = f(x_k) - \alpha(1-\frac{\alpha L B^2}{2})||f'(x_k)||^2.
\end{aligned}
\end{equation}
This inequality shows that if $x_k$ is not a minimizer, then the stochastic gradient descent iteration is expected to decrease the objective function for any step size satisfying
\begin{equation}
\label{eq:eleven}
0 < \alpha < \frac{2}{LB^2}.
\end{equation}

\section{Linear Convergence for Strongly Convex Objectives}

We now use the bound~\eqref{eq:ten} to show that, for strongly convex functions, constant step sizes satisfying~\eqref{eq:eleven} lead to an expected linear convergence rate.  First, use $x = x_k$ in~\eqref{eq:twelve} and minimize both sides of~\eqref{eq:twelve} with respect to $y$ to obtain
\[
f(x^*) \geq f(x_k) - \frac{1}{2\mu}||f'(x_k)||^2,
\]
where $x^*$ is the minimizer of $f$. Subsequently, we have
\[
-||f'(x_k)||^2 \leq -2\mu(f(x_k)-f(x^*)).
\]
Now use this in~\eqref{eq:ten} and assume the step sizes satisfy~\eqref{eq:eleven} to get
\[
\mathbb{E}[f(x_{k+1})] \leq f(x_k) - 2\mu\alpha(1 - \frac{\alpha L B^2}{2})[f(x_k) - f(x^*)].
\]
We now subtract $f(x^*)$ from both sides and take the expectation with respect to the sequence $\{e_0,e_1,\dots,e_{k-1}\}$ to obtain
\begin{align*}
\mathbb{E}[f(x_{k+1})] - f(x^*) &\leq \mathbb{E}[f(x_k)] - f(x^*) - 2\mu\alpha(1 - \frac{\alpha L B^2}{2})[\mathbb{E}[f(x_k)] - f(x^*)]\\
& = \left(1 - 2\mu\alpha(1 - \frac{\alpha L B^2}{2})\right)[\mathbb{E}[f(x_k)] - f(x^*)].
\end{align*}
Applying this recursively we have
\[
\mathbb{E}[f(x_k)] - f(x^*) \leq \rho^k[f(x_0)-f(x^*)],
\]
for some $\rho < 1$. Thus, the difference between the expected function value $\mathbb{E}[f(x_k)]$ and the optimal function value $f(x^*)$ decreases \emph{geometrically} in the iteration number $k$.

In the particular case of $\alpha = \frac{1}{LB^2}$, this expression simplifies to
\[
\mathbb{E}[f(x_k)] - f(x^*) \leq \left(1 - \frac{\mu}{LB^2}\right)^k[f(x_0)-f(x^*)],
\]
and thus the method approaches the $(1-\mu/L)^k$ rate of the deterministic method with a step size of $1/L$~\citep[see][\S8.6]{luenberger2008linear} as $B$ approaches one.

\section{Sublinear $O(1/k)$ Convergence for Convex Objectives}

We now turn to the case where $f$ is convex but not necessarily strongly convex.  In this case, we show that if at least one minimizer $x^*$ exists, then a step size of $\alpha = \frac{1}{LB^2}$ leads to an $O(1/k)$ error.  By convexity, we have for any minimizer $x^*$ that
\[
f(x_k) \leq f(x^*) + \langle f'(x_k) , x_k - x^* \rangle,
\]
and thus for any $\beta \leq 1$ that 
\[
f(x_k) \leq \beta f(x_k) + (1-\beta)f(x^*) + (1-\beta)\langle f'(x_k) , x_k - x^* \rangle.
\]
We use this to bound $f(x_k)$ in~\eqref{eq:13} to get
\begin{equation}
\label{eq:14}
\begin{aligned}
f(x_{k+1}) & \leq \beta f(x_k) + (1-\beta)f(x^*) + (1-\beta)\langle f'(x_k) , x_k - x^* \rangle \\
&\hspace*{1cm} - \alpha(1 - \frac{\alpha L}{2})||f'(x_k)||^2 - \alpha(1 - \alpha L)\langle f'(x_k),e_k\rangle + \frac{\alpha^2L}{2}||e_k||^2.
\end{aligned}
\end{equation}
Note that 
\begin{align*}
\frac{1}{2\alpha} \left(\|x_k - x^*\|^2 - \|x_{k+1} - x^*\|^2\right) & =
\frac{1}{2\alpha} \left(\|x_k - x^*\|^2 - \|x_k - \alpha f'(x_k) - \alpha e_k - x^*\|^2\right) \\
& = -\frac{\alpha}{2} \|f'(x_k)\|^2 -\frac{\alpha}{2} \|e_k\|^2 - \alpha \langle f'(x_k), e_k \rangle\\
	&\hspace*{1cm} + \langle f'(x_k), x_k - x^* \rangle + \langle e_k, x_k - x^* \rangle \; ,
\end{align*}
and using this to replace $\langle f'(x_k), x_k - x^* \rangle$ in~\eqref{eq:14} we obtain the ugly expression
\begin{align*}
f(x_{k+1}) & \leq \beta f(x_k) + (1 - \beta) f(x^*) + \frac{1 - \beta}{2\alpha}\left(\|x_k - x^*\|^2 - \|x_{k+1} - x^*\|^2\right)\\
&\hspace*{1cm} +\frac{\alpha(1 - \beta)}{2}\left(\|f'(x_k)\|^2 + \|e_k\|^2\right) + (1 - \beta)\alpha \langle f'(x_k), e_k \rangle - (1 - \beta) \langle e_k, x_k - x^*\rangle\\
&\hspace*{1cm} - \alpha(1 - \frac{\alpha L}{2})||f'(x_k)||^2 - \alpha (1 - L\alpha) \langle f'(x_k), e_k \rangle + \frac{L\alpha^2}{2}\|e_k\|^2 \; .
\end{align*}
Taking the expectation with respect to $e_k$ and using properties~\eqref{eq:eight} and~\eqref{eq:nine}, this becomes
\begin{equation}
\label{eq:15}
\begin{aligned}
E[f(x_{k+1})] & \leq \beta f(x_k) + (1 - \beta) f(x^*) + \frac{1 - \beta}{2\alpha}\left(\|x_k - x^*\|^2 - E[\|x_{k+1} - x^*\|^2]\right)\\
&\hspace*{1cm} +\frac{\alpha(1 - \beta)}{2}\left(\|f'(x_k)\|^2 + (B^2-1)\|f'(x_k)\|^2\right)\\
&\hspace*{1cm} - \alpha(\frac{2-\alpha L}{2})||f'(x_k)||^2 + \frac{L\alpha^2(B^2-1)}{2}\|f'(x_k)\|^2 \; .
\end{aligned}
\end{equation}
Using $\alpha = \frac{1}{LB^2}$, we can make all terms in $\|f'(x_k)\|$ cancel out by choosing $\beta = 1 - \frac{1}{B^2}$ because
\[
\alpha(1 - \beta) B^2 - 2\alpha +L\alpha^2B^2 = \alpha - 2\alpha + \alpha = 0.
\]
We now take the expectation of~\eqref{eq:15} with respect to $\{e_0,e_1,\dots,e_{k-1}\}$ and note that $(1-\beta)/\alpha = L$ to obtain
\[
E[f(x_{k+1})]- f(x^*) \leq \beta E[f(x_k)] - \beta f(x^*) + \frac{L}{2}\left(E[\|x_k - x^*\|^2] - E[\|x_{k+1} - x^*\|^2]\right).
\]
If we sum up the error from $k =0$ to $(n - 1)$, we have
\begin{align*}
\sum_{k=0}^{n-1} \left(E[f(x_{k+1})] - f(x^*)\right) & \leq \beta \sum_{k=0}^{n-1} \left(E[f(x_k)] - f(x^*)\right) + \frac{L}{2} \left(\|x_0 - x^*\|^2 - E[\|x_n - x^*\|^2]\right)\\
	& \leq \beta \sum_{k=1}^{n} \left(E[f(x_{k})] - f(x^*)\right) + \beta\left(f(0) - f(x^*)\right) + \frac{L}{2} \|x_0 - x^*\|^2 \; .
\end{align*}
Hence, we have
\begin{align*}
(1 - \beta) \sum_{k=0}^{n-1} \left(E[f(x_{k+1})] - f(x^*)\right) & \leq \beta\left(f(0) - f(x^*)\right) + \frac{L}{2} \|x_0 - x^*\|^2 \; .
\end{align*}
Since $E[f(x_{k+1})]$ is a non-increasing function of $k$, the sum on the left-hand side is larger than $k$ times its last element. Hence, we get
\begin{align*}
E[f(x_{k+1})] - f(x^*) & \leq \frac{1}{k}\sum_{i=0}^{k-1}(E[f(x_{i+1})] - f(x^*))\\
 & \leq \frac{\beta\left(f(0) - f(x^*)\right) + \frac{L}{2} \|x_0 - x^*\|^2}{k(1 - \beta)}\\
 & = \frac{2(B^2-1)\left(f(0) - f(x^*)\right) + LB^2 \|x_0 - x^*\|^2}{2k}\\
 & = O(1/k).
\end{align*}

\bibliographystyle{abbrvnat}
\bibliography{bib}

\end{document}